\documentclass[12pt]{article}
\usepackage{amsmath,amsfonts,amssymb,amsthm}
 \usepackage{pst-plot}  
\newcommand{\I}{{\bf 1}}
\newtheorem{proposition}{Proposition}[section]
\newtheorem{theorem}[proposition]{Theorem}
\newtheorem{corollary}[proposition]{Corollary}
\newtheorem{lemma}[proposition]{Lemma}
\newtheorem{remark}[proposition]{Remark}

\newtheorem{assumption}[proposition]{Assumption}
\newtheorem{scenario}{Scenario}

\newcommand{\nc}{\newcommand}
\nc{\R}{{\mathbb R}}
\nc{\N}{{\mathbb N}}
\nc{\Z}{{\mathbb Z}}
\DeclareMathOperator{\card}{card}

\DeclareMathOperator{\sgn}{sgn}
\nc{\BP}{\mathbb{P}}
\nc{\BE}{\mathbb{E}}
\nc{\BQ}{\mathbb{Q}}
\nc{\bN}{{\mathbf N}}

\newcommand{\dle}{\stackrel{d}{\le}}             
\newcommand{\deq}{\stackrel{d}{=}}               
\newcommand{\dto}{\xrightarrow{d}}               
\newcommand{\oxi}{\overline\xi}                  
\newcommand{\uxi}{\underline\xi}                 
\newcommand{\omu}{\overline\mu}                  
\newcommand{\olambda}{\overline\lambda}          
\newcommand{\oS}{\overline S}                    
\newcommand{\uS}{\underline S}                   
\newcommand{\umu}{\underline\mu}                 
\newcommand{\ulambda}{\underline\lambda}         

\begin{document}

\author{K.A. Borovkov\footnote{
Department of Mathematics and Statistics, University of Melbourne,
k.borovkov@ms.unimelb.edu.au}
\ and \ G. Last
\footnote{
Institut f\"ur Stochastik, Universit\"at Karlsruhe (TH),
last@math.uni-karlsruhe.de} }
\title{On level crossings
for a general class of\\ piecewise-deterministic Markov processes}
\date{}
\maketitle
\begin{abstract}
\noindent
We consider a piecewise-deterministic Markov process  
governed by a jump intensity function, a rate function that determines
the behaviour between jumps, and a stochastic kernel describing
the conditional distribution of jump sizes.
We study the point process $N^b_+$ of upcrossings 
of some level $b$. Our main result shows that
a suitably scaled point process $N^b_+(\nu(b)t)$, $t\ge 0$,
converges, as $b\to\infty$, weakly to a geometrically compound Poisson process.
We also prove a version of Rice's formula
relating the stationary density of the process to level crossing
intensities. This formula provides an interpretation
of the scaling factor $\nu(b)$. While our proof of the limit theorem requires
additional assumptions, Rice's formula holds whenever
the (stationary) overall intensity of jumps is finite.

\end{abstract}

\vspace{0.1cm}

\section{Introduction}
\setcounter{equation}{0}

We consider a real-valued piecewise-deterministic
Markov process $(X_t)_{t\ge 0}$ whose distribution is determined by 
a drift coefficient $\mu:\R\rightarrow\R$, a jump intensity function 
$\lambda:\R\rightarrow [0,\infty)$, and a stochastic kernel $J(x,dz)$ from 
$\R$ to $\R$. 
The process $(X_t)$ is right-continuous and jumps at (positive) epochs 
$T_1<T_2<\ldots$. Between the jumps it moves along an integral curve
determined by $\mu$. We assume that $\mu$ is right-continuous
and that $D_\mu:=\{u:\mu(u)=0\}$ is a locally finite set. 
The occurence of jumps is governed by the stochastic jump intensity
$\lambda(X_t)$. Given the $n$-th jump epoch $T_n$, the conditional
distribution of the size $Z_n$ of the $n$-th jump 
is $J(X_{T_n-},\cdot)$, where $X_{t-}$ is the 
value of the process just before $t>0$.
We will assume that the process has an invariant distribution $\pi$
and refer to the Appendix for conditions guaranteeing the existence of
a unique stationary distribution. It is then
essentially well-known (\cite{Costa}, \cite{Zheng})
that the stationary distribution $\pi$ is absolutely continuous
on $\R\setminus D_\mu$, and we let $p$ denote its density.
We note that $\pi$ might have atoms in $D_\mu$.

The process $(X_t)$ is a generic model of applied probablity. Special cases
have been extensively studied in the literature. We just mention
storage processes (\cite{HR},\cite{MT1993}),    
stress release models (\cite{BorovkovV}, \cite{VereJones1988},\cite{Zheng}),
queueing models (\cite{BrowneSigman},\cite{MT1993}), and
repairable systems (\cite{LS98a}). It is mostly assumed that
$J(x,\cdot)$ does not depend on $x\in\R$ and that the jumps
are either only non-negative or only non-positive. An extensive discussion
of several  ergodicity properties for a constant (positive) $\mu$ and
negative jumps is given in \cite{La04}.
General properties of piecewise-deterministic Markov processes are studied in
\cite{Davis}.

Now assume that $X_0$ has the distribution $\pi$.
Then $(X_t)$ is a {\em stationary process}, and  
the sequence $(T_n)$ forms a {\em stationary point process}. We assume
that the {\em intensity} of $(T_n)$ (the expected number of points in an interval of
unit length) is finite. Again we refer to the Appendix for an explicit
assumption that is sufficient for this finiteness. 
We say that $(X_t)$ has an {\em upcrossing}
(resp.\ {\em downcrossing})
of level $u\in\R$ at time $s>0$ if there is some
$\delta>0$ such that $X_t<u$ (resp.\ $X_t\ge u$) for $s-\delta\le t<s$ and
$X_t\ge u$ (resp.\ $X_t< u$) for $s< t\le s+\delta$. 
If, in addition, $X_{s-}=X_s(=u)$ then we speak of a {\em continuous upcrossing}
(resp.\  {\em downcrossing}). It is easy to see that the
set of all continuous up- and downcrossings forms a stationary  
point process $N^u$.  Note that there are no continuous downcrossings in
case $\mu(u)> 0$ and no continuous upcrossings in case $\mu(u)<0$. 
The intensity of $N^u$ is denoted by $\nu(u)$. 
As the intensity of $(T_n)$ is assumed finite, it is easy to see
that $\nu(u)$ is finite for any $u\in\R$.

Our first aim in this paper is to prove the following version of {\em Rice's formula}:
\begin{align}\label{rice}
\nu(u)=|\mu(u)|p(u),\quad u\notin D_\mu. 
\end{align}
The simplicity of this formula is striking.  If $(X_t)$ is
ergodic, \eqref{rice} can be explained by looking at the long-run proportion
of time that $(X_t)$ spends in an infinitesimal interval containing $u$.
Formula \eqref{rice} is a direct analog of the classical Rice formula \cite{Rice},
which holds for smooth processes
and plays a rather important role in engineering.
A rigorous treatment of Rice's formula is given in \cite{Lead66} 
and a more recent discussion in \cite{LeadSpan02}.
An analog of \eqref{rice} for (discontinuous) Poisson shot-noise processes
has been studied in \cite{BN72}.

Let $\nu_{+,d}(b)$ and $\nu_{-,d}(b)$ denote resp.\ the intensities of {\em discontinuous}
up- and downcrossings of the level $b$. Our proof of  \eqref{rice} 
uses the simple relation $\nu(u)=|\nu_{+,d}(u)-\nu_{-,d}(u)|$, see Lemma \ref{l1}.
In fact, \eqref{rice} can be rewritten as
\begin{align}\label{rice2}
\nu_{-,d}(u)-\nu_{+,d}(u)=\mu(u)p(u),\quad u\notin D_\mu. 
\end{align}
Such equalities for level
crossing intensities are widely used in queueing theory.
We refer here to the early reference \cite{BrillPosner}
and the survey \cite{Doshi92}. 
It is quite remarkable that the queueing literature does not
take notice of the close relationship between \eqref{rice2}
and the results in \cite{Rice} (or \cite{BN72}).
Equation \eqref{rice2} is mostly derived for Poisson driven models. 
In principle, the level crossing method can also be applied in
more general cases (see e.g.\ \cite{Doshi92}). There are, however,
many implicit model assumptions, that make a direct derivation
of \eqref{rice} non-trivial. 
So to the best of our knowledge, the result \eqref{rice} must be
considered as new. Moreover, we will establish this formula 
under a minimal set of assumptions. In particular, the existence
of the stationary density need not be assumed, but is a consequence of our
model assumptions. Even ergodicity is not needed.

Our second and main aim in this paper is to derive limit results for the
point process $N^b_+$ of all upcrossings of the level $b\to\infty$.
Whenever the intensity $\nu_+(b)$ of $N^b_+$ is positive, 
we introduce the scaled point process 
$M^b(t):=N^b_+(\nu_+(b)^{-1}t)$, $t\ge 0$.
It is stationary and has intensity $1$. Under our assumptions
(see the scenarios below), equation \eqref{rice} will imply that the intensity $\nu_+(b)$
can be explicitly expressed as
\begin{align}\label{rb}
\nu_+(b)=|\mu(b)|p(b),
\end{align}
for all sufficiently large $b$. 
We refer here to Section \ref{aslevel} for more details.
We will study the limiting behavior of $M^b$ under the following
three scenarios and some additional assumptions, see
\eqref{NegDrift1}-\eqref{NegDrift3}.

\begin{scenario}\label{s1}\rm
We have $\mu(y)\to -\infty$ as
$y\to\infty$, and there exists a $u_0\in\R$ such that $J(x, (-\infty,0))=0$ for
$x\ge u_0$ (no negative jumps from states $x\ge u_0$).
\end{scenario}

\begin{scenario}\label{s2}\rm
We have $\lambda(y)\to  \infty$ as
$y\to\infty$, $\mu(y)$ is positive for all sufficiently large $y$, and 
$J(x,(0,\infty))=0$ for all $x\in \R$ (no positive jumps). 
\end{scenario}

\begin{scenario}\label{s3}\rm
As $y\to\infty$ we have $\mu(y)\to\mu (\infty)\in \R\setminus\{0\}$ and  
$\lambda(y)\to \lambda(\infty)\in [0,\infty)$.
In case $\mu(\infty)<0$ there exists a $u_0\in\R$ such that $J(x, (-\infty,0))=0$ for
$x\ge u_0$ and in case $\mu(\infty)>0$ we have $J(x,(0,\infty))=0$ for all $x\in\R$.
Moreover, $J(y,\cdot)$ converges weakly, as $y\to\infty$, to a probability measure
$J(\infty,\cdot)$ on $\R$.
\end{scenario}

In the first two scenarios the point process $M^b$ will converge, as $b\to\infty$,
in distribution to a Poisson process. The explanation of
this phenomenon is quite simple. Fixing a level $u>u_0$,
the trajectory of the process $(X_t)$ can be split into i.i.d.\ cycles
between the successive continuous crossings of this level.
Then hitting a high level $b$ during a particular cycle will be a `rare event'. 
In both scenarios, with a probability arbitrary close to 1 for large enough
$b$, given that the level $b$ was exceeded during a cycle there is exactly one  
upcrossing of that level during this cycle.

In the third scenario the limiting behaviour of $M^b$ is slightly more complicated.
The crossing of a high level $b$ is still a rare event. However, 
given the level $b$ was exceeded during a  cycle, the conditional
distribution of the number of continuous crossings of that level during this 
cycle will be geometric with a parameter that
converges as $b-u\to\infty$ to some number $\rho\in(0,1)$.
Therefore the limit is a {\em geometrically compound Poisson process} $\Pi_\rho$
which is defined as follows. Each point of a homogeneous Poisson process of
intensity $(1-\rho)$ gets (independently of the other points)
a mass $k\in\{1,2,\ldots\}$ with probability $(1-\rho)\rho^{k-1}$. The 
resulting stationary point process $\Pi_\rho$  has {\em independent increments}
and geometrically distributed multiplicities.
As the above geometric distribution has mean $1/(1-\rho)$, the intensity
of $\Pi_\rho$ is $1$. 

For Gaussian processes it is well-known that the point process
of times of crossing a high level 
are asymptotically Poisson, see e.g.\  \cite{LinLeadRoo83} and the references
given there. But to the best of our knowledge the present paper is the
first to establish such a limit theorem for jump processes.
Compound Poisson limits for exceedances and upcrossings of sequences 
are summarized in \cite{FHR04}.
A general compound Poisson limit theorem for strongly mixing random measures
has been derived in \cite{LeadHsing90}. We are not aware of any straightforward
way to derive our theorem from these results.
An immediate consequence of our main theorem is that the first time of
crossing a high level is asymptotically exponentially distributed.
A discussion of this well-known phenomenon can be found, for instance,
in \cite{Aldous89} and Section VI.4 of \cite{As03}.
However, in the present framework the result seems to be new.

This paper is organized as follows. Section 2 contains the detailed
definition of the process as well as some of its fundamental
properties. Section 3 provides the proof of Rice's formula.
The Poisson limit theorem is the topic of the final and main Section 4.

\section{Definition and basic properties of the process}\label{sec2}
\setcounter{equation}{0}

We consider a right-continuous function $\mu:\R\rightarrow\R$ 
such that the set $D_\mu$ of zeros of $\mu$ is locally finite. 
We assume that, for any $x\in\R$, there
exists a unique continuous function  $q(x,\cdot):[0,\infty)\rightarrow\R$ 
satisfying the integral equation
\begin{align}\label{ie}
q(x,t)=x+\int^t_0 \mu(q(x,s))\, ds,\quad t\ge 0.
\end{align}
The jump intensity $\lambda$ is asumed to be measurable, locally bounded
and such that
\begin{align}\label{infty}
\int^\infty_0\lambda(q(x,s))ds=\infty,\quad x\in\R.
\end{align}
For the jump distribution we assume that $J(x,\{0\})=0$ for all 
$x\in\R$ (see also Remark \ref{r1}).

Formally, our process $(X_t)$ is defined as follows. We consider a 
measurable space $(\Omega,\mathcal{F})$ that is rich enough to carry
a {\em marked point process} $\Phi=((T_n,Z_n))_{n\ge 1}$ on $[0,\infty)$ with real-valued
random variables (marks) $Z_n$ and a real-valued random variable $X_0$. Between
the jumps the process is defined by $X_t:=q(X_0,t)$
on $[0,T_1)$ and $X_t=q(X_{T_n},t-T_n)$ on $[T_n,T_{n+1})$, 
$n\ge 1$. At the jump epochs $T_n$ we have $X_{T_n}:=X_{T_{n}-}+Z_n$,
where $X_{T_{n}-}:=\lim_{s\to T_{n}-}X_s=q(X_{T_{n-1}},T_n-T_{n-1})$.
Finally, we define $X(t):=\Delta$ for $t\ge T_\infty$,
where $\Delta$ is a point external to $\R$ and $T_\infty:=\lim_{n\to\infty}T_n$.

For any probability measure $\sigma$ on $\R$ we consider a probability measure
$\BP_\sigma$ on $(\Omega,\mathcal{F})$ such that $\BP_\sigma(X_0\in\cdot)=\sigma$
and the following  properties hold.
The conditional distribution of $T_1$ given $X_0$ is specified by
\begin{align}\label{T1}
\BP_\sigma(T_1\le t|X_0)=1-\exp\left[-\int^t_0\lambda(q(X_0,s))ds\right]\quad 
\BP_\sigma\text{-a.s.}
\end{align}
Similarly we assume for $n\ge 1$  that, $\BP_\sigma$-almost surely, 
\begin{align}\label{Tn}
\BP_\sigma(T_{n+1}-T_n\le t|X_0,T_1,Z_1,\ldots,T_n,Z_n)
=1-\exp\left[-\int^t_0\lambda(q(X_{T_n},s))ds\right].
\end{align}
By \eqref{infty} the jump epochs $T_n$ are indeed all finite a.s.
The conditional distributions of the jump sizes are given by
\begin{align}\label{Zn}
\BP_\sigma(Z_{n+1}\in\cdot |X_0,T_1,Z_1,\ldots,T_n,Z_n,T_{n+1})
=J(X_{T_{n+1}-},\cdot) \quad \BP_\sigma\text{-a.s.},\, n\ge 0.
\end{align}
Since $J(x,\{0\})=0$, $x\in\R$, we can assume that $Z_n(\omega)\ne 0$ for all $n\ge 1$
and $\omega\in\Omega$.

The conditional distribution of $\Phi$ given $X_0$ is now completely specified.
Our assumptions imply that $(X_t)$ is a homogeneous Markov process with respect to
the family $\{\BP_x:x\in\R\}$, where $\BP_x:=\BP_{\delta_x}$ is the measure
correponding to the initial distribution supported by $x$. 
The expectations with respect to $\BP_\sigma$ and $\BP_x$ are denoted by
$\BE_\sigma$ and $\BE_x$ respectively. 
Actually, $(X_t)$ is piecewise-deterministic Markov process in the
terminology of \cite{Davis}.

\begin{remark}\label{r1}\rm As $Z_n\ne 0$ for all $n\ge 1$ there is a one-to-one
correspondence between $\Phi$ and $(X_t)$. The former condition can be easily
dispensed with by suitably augmenting the process $(X_t)$. 
\end{remark}

\begin{remark}\label{r0}\rm In many applications (queueing and dam models, repairable
systems) the process $(X_t)$ is non-negative, in the
sense that $X_t\ge 0$ for all $t\ge 0$ whenever $X_0\ge 0$. Such a situation can be
accomodated by choosing the characteristics so that $(-\infty,0)$
becomes {\em transient} for the process. A possible choice is $\mu(x)=1$ and
$\lambda(x)=0$ for $x<0$. Any stationary distribution of $(X_t)$ is then
concentrated on $[0,\infty)$.
\end{remark}

\begin{remark}\label{r12}\rm We assumed that the solution $q(x,t)$ 
of \eqref{ie} is defined for all $t\ge 0$. This could be generalized as follows.
Suppose that for any $x\in\R$ there is a $t_\infty(x)\in(0,\infty]$ such
that $q(x,\cdot)$ is the unique continuous function on $[0,t_\infty(x))$  
satisfying \eqref{ie} for all $t\in [0,t_\infty(x))$. Assuming instead of
\eqref{infty} that
$\int^{t_\infty(x)}_0\lambda(q(x,s))ds=\infty$, $x\in\R$,
we can still use \eqref{T1}, \eqref{Tn}, and \eqref{Zn} to define a marked point process
$\Phi$ such that a.s.\ $T_1<t_\infty(X_0)$ and
$T_{n+1}-T_n<t_\infty(X_{T_n})$, $n\ge 1$.
Hence we can define the Markov process $(X_t)$ as before.
All results of this paper remain valid in this more general framework.
\end{remark}

The next result provides the (generalized) {\em infinitesimal generator} of $(X_t)$.
Set
$$
\tau_m:=\inf\{t\ge 0:|X(t)|\ge m\},\quad m\in\N.
$$

\begin{proposition}\label{p0}
Let $f:\R\rightarrow \R$ be absolutely continuous    
with a Radon--Nikodym derivative $f'$ and let $f'$ as well as 
the function $x\mapsto\lambda(x)\int (f(x+z)-f(x))J(x,dz)$ be locally bounded. Then,
for any probability measure $\sigma$ on $\R$,
\begin{multline}\label{13}
\BE_\sigma f(X_{t\wedge \tau_m})=\BE_\sigma f(X_0)+\BE_\sigma\int^{t\wedge\tau_m}_0 f'(X_s)\mu(X_s)ds\\
+\BE_\sigma\int_0^{t\wedge\tau_m}\int_\R (f(X_s+z)-f(X_s))\lambda(X_s)J(X_s,dz)ds.
\end{multline}
\end{proposition}
{\sc Proof:} Denote by $(\mathcal{F}_t)$ the filtration  generated by
$X_0$ and the restriction of $\Phi$ to $[0,t]\times \R$. 
Using basic results on marked point processes (see e.g.\ chapter 4 in \cite{LastBrandt})
we obtain from \eqref{T1}, \eqref{Tn}, and \eqref{Zn} that
\begin{align}\label{comp}
\BE_\sigma\sum_{n=1}^\infty h(T_n,Z_n)=
\BE_\sigma\int^\infty_0\int_\R h(t,z) \lambda(X_t)J(X_t,dz)dt,
\end{align}
for all {\em predictable} 
$h:\Omega\times [0,\infty)\times[0,\infty)\rightarrow [0,\infty)$.
We can now proceed as  in Section 8 of \cite{La04} to obtain the result.\qed

\vspace{0.3cm}
We have to make two basic assumptions on the process.
They will be discussed in the Appendix.

\begin{assumption}\label{a1}\rm We have $\BP_x(T_\infty=\infty)=1$ for all
$x\in\R$, and the process $(X_t)$ has an invariant distribution $\pi$.
\end{assumption}

In view of Remark \ref{r1}, the
marked point process $\Phi$ is stationary under $\BP_\pi$, see
\cite{BB} for more detail on this stationarity.
In particular, the distribution of  $(N(t+s)-N(s))_{t\ge 0}$ does not depend on
$s\ge 0$, where $N(t):=\card\{n\ge 1:T_n\le t\}$ is the number of jumps
in the time interval $[0,t]$. The (stationary) {\em intensity} of $N$
is defined by 
\begin{align*}
\lambda_\pi:=\BE_\pi N(1).
\end{align*}
 

\begin{assumption}\label{a2}\rm We have $\lambda_\pi<\infty$.
\end{assumption}

Let $g:[0,\infty)\times\R\times\R\rightarrow [0,\infty)$ be measurable.
Applying \eqref{comp} with $\sigma=\pi$ and $h(t,z):=g(t,X_{t-},z)$, and using Fubini's 
theorem, we obtain
\begin{align}\label{94}
\BE_\pi\sum^\infty_{n=1}g(T_n,X_{T_n-},Z_n)
=\int^\infty_0\int_\R\int_\R g(s,x,z)\lambda(x)J(x,dz)\pi(dx)ds.
\end{align} 
Choosing $g(s,x,z)=\I\{0\le s\le 1\}$ we obtain the equality in
\begin{align}\label{intensity}
\lambda_\pi=\int\lambda(x)\pi(dx)<\infty.
\end{align}

A quick consequence of Proposition \ref{p0}  is the
following  (basically well-known) integral equation for $\pi$.

\begin{proposition}\label{p1} 
Let $f:\R\rightarrow \R$ be bounded and absolutely continuous    
with a continuos Radon--Nikodym derivative $f'$ that has a compact support. 
Then
$$
\int f'(x)\mu(x)\pi(dx)=\iint\lambda(x)(f(x)-f(x+z))J(x,dz)\pi(dx).
$$ 
\end{proposition}
{\sc Proof:} The assumptions on $f$ allow to use formula \eqref{13}. 
Because of Assumption \ref{a1} the process $(X_t)$ is real-valued
and locally bounded. Hence we have $\BP_\pi$-a.s.\ that
$\tau_m\to\infty$ as $m\to\infty$. As $f$ is
bounded, the left-hand side of \eqref{13} converges to 
$\BE_\pi f(X_{t})=\BE_\pi f(X_{0})$. 
As $f'$ has a compact support and $\mu$ is locally bounded, the second
term on the right-hand side of \eqref{13} converges as well.
For the third term we can use \eqref{intensity} and bounded convergence 
to conclude that
$$
0=\BE_\pi\int^{t}_0 f'(X_s)\mu(X_s)ds
+\BE_\pi\int_0^{t}\int_\R \lambda(X_s)(f(X_s+z)-f(X_s))J(X_s,dz)ds.
$$
Using Fubini's theorem and stationarity again, we obtain the assertion.\qed
\vspace{0.3cm}

Some relationships between
$\pi$ and the stationary distribution of the imbedded process $(X_{T_n})$
can be found in \cite{Costa}. 

\section{Rice's formula}\label{sproof}
\setcounter{equation}{0}

In this section we will prove the following assertion, establishing the Rice 
formula \eqref{rice}:

\begin{theorem}\label{trice} Under Assumptions \ref{a1} and \ref{a2} the 
stationary distribution $\pi$ has a right-continuous density $p$ on $\R\setminus D_\mu$
satisfying $\nu(u)=|\mu(u)|p(u)$ for all $u\notin D_\mu$.
\end{theorem}

We prepare the proof with an auxiliary result and start with
introducing some notation.
We say that $(X_t)$ has a {\em discontinuous upcrossing} (resp.\
{\em discontinuous downcrossing})
of level $u$ at time $s>0$ if $X_{s}\ge u> X_{s-}$ 
(resp.\ $X_{s-}\ge u> X_s$).
The point processes of these discontinuous down- and upcrossings
are denoted by $N^u_{+,d}$ and $N^u_{-,d}$. 
In this section we take $\BP_\pi$ to be the underlying probability measure.
Then $\Phi$ is a stationary marked point
process, and $N^u_{+,d}$ and $N^u_{-,d}$ are (jointly) stationary point processes. 
Their intensities are denoted by $\nu_{+,d}(u)$ and
$\nu_{-,d}(u)$, respectively.

\begin{lemma}\label{l1} For any  $u\in\R$ we have
$\nu(u)=\nu_{+,d}(u)-\nu_{-,d}(u)$ in case $\mu<0$ on $(u,u+\varepsilon)$
for some $\varepsilon>0$
and $\nu(u)=\nu_{-,d}(u)-\nu_{+,d}(u)$ in case $\mu> 0$
on $(u,u+\varepsilon)$ for some $\varepsilon>0$.
\end{lemma}
{\sc Proof:} Assume that $\mu<0$ on $(u,u+\varepsilon)$ for some $\varepsilon>0$. 
(The argument for the other case is the same.)
As the solution of \eqref{ie} is unique and $\mu$ is right-continuous,
there are no continuous upcrossings of level $u$. Therefore,
between any two (discontinuous or
continuous) successive downcrossings there must be exactly one 
discontinuous upcrossing of $u$. 
Hence we have for any $t\ge 0$ that
$$
N^u_{+,d}(t)-1\le N^u_{-,d}(t)+N^u(t)\le N^u_{+,d}(t)+1.
$$  
Taking expectations gives
$$
\nu_{+,d}(u)t-1\le \nu_{-,d}(u)t+\nu(u)t\le \nu_{+,d}(u)t+1.
$$
Dividing by $t$ and letting $t\to\infty$, yields the assertion. \qed

\vspace{0.3cm}
{\sc Proof of Theorem \ref{trice}:} Let $u\in \R$. 
Choosing  $g(s,x,z):=\I\{0\le s\le 1,x<u\le x+z\}$
(resp.\ $g(s,x,z):=\I\{0\le s\le 1,x\ge u> x+z\}$) in \eqref{94} yields 
\begin{align}\label{nu+}
\nu_{+,d}(u)&=\iint \I\{x<u\le x+z\}\lambda(x)J(x,dz)\pi(dx),\\
\label{nu-}
\nu_{-,d}(u)&=\iint \I\{x\ge u> x+z\}\lambda(x)J(x,dz)\pi(dx).
\end{align}

Let $f$ be a function satisfying the assumptions of Proposition \ref{p1}.
By \eqref{nu+} and \eqref{nu-} we have
\begin{align*}
\int f'(u)&(\nu_{-,d}(u)-\nu_{+,d}(u))du=
\iiint f'(u) \I\{x>u\ge x+z\}du\lambda(x)J(x,dz)\pi(dx)\\
&-\iiint f'(u) \I\{x+z>u\ge x\}du\lambda(x)J(x,dz)\pi(dx)\\
=&\iint \I\{z<0\}(f(x)-f(x+z))\lambda(x)J(x,dz)\pi(dx)\\
&-\iint \I\{z>0\}(f(x+z)-f(x))\lambda(x)J(x,dz)\pi(dx)\\
=&\iint (f(x)-f(x+z))\lambda(x)J(x,dz)\pi(dx).
\end{align*}
Therefore we obtain from Proposition \ref{p1} that
\begin{align*}
\int f'(u)(\nu_{-,d}(u)-\nu_{+,d}(u))du=
\int f'(u)\mu(u)\pi(du).
\end{align*}
The class of functions $f'$ that are allowed in the above formula is rich 
enough to conclude first that $\pi$ is absolutely continuous on $\R\setminus D_\mu$
and second, that the density $p$ satisfies
\begin{align}\label{88}
\mu(u)p(u)=\nu_{-,d}(u)-\nu_{+,d}(u)
\end{align}
for almost all $u\notin D_\mu$. 
By \eqref{nu+} and \eqref{nu-} the function $\nu_{-,d}-\nu_{+,d}$ is  left-continuous
so that Lemma \ref{l1} shows that $\nu$ is left-continuous on $\R\setminus D_\mu$.
In fact, for $u\notin D_\mu$ the lemma remains true, 
if  $\nu_{-,d}$ and $\nu_{+,d}$ are replaced
by the corresponding right-continuous versions. Hence $\nu$ is even continuous
on $\R\setminus D_\mu$, and we can use \eqref{88} to redefine a right-continuous density $p$. 
Lemma \ref{l1} implies the assertion.
\qed

\section{Asymptotics of level crossings}\label{aslevel}
\setcounter{equation}{0}

In this section we write $\BP:=\BP_\pi$. 
Consider the point process $N^b_+$ of all
upcrossings of some level $b\in\R$ and let $\nu_+(b)$ denote its intensity
(under $\BP$). It is given by 
$$
\nu_+(b)=\I\{\mu(b)> 0\}\nu(b)+\nu_{+,d}(b),\quad b\notin D_\mu, 
$$
where we refer to the Introduction and Section \ref{sproof} for the definition
of the intensities $\nu(b)$ and $\nu_{+,d}(b)$. 
>From Lemma \ref{l1}, \eqref{nu+}, and  \eqref{nu-} we obtain that $\nu_+(b)\to 0$ as $b\to\infty$.
If $\mu(x)<0$ and $J(x,(-\infty,0))=0$ for all
$x\ge u_0$ (no negative jumps from levels above $u_0$),
we conclude from \eqref{nu-} and Lemma \ref{l1} that $\nu_{-,d}(b)=0$ and 
$\nu_+(b)=\nu_{+,d}(b)=\nu(b)$ for $b\ge u_0$.
If $J(x,[(u_0-x)^+,\infty))=0$ for all $x\in\R$  (no positive jumps to levels
above $u_0$) and $\mu(x)> 0$ for all $x\ge u_0$,
we conclude from \eqref{nu+} and Lemma \ref{l1} that $\nu_{+,d}(b)=0$ and 
$\nu_+(b)=\nu_{-,d}(b)=\nu(b)$ for $b\ge u_0$. (Here $a^+:=\max\{a,0\}$ denotes the
positive part of $a$.)
In either case, Theorem \ref{trice} implies that
\eqref{rb} holds for $b\ge u_0$. 
Whenever $\nu_+(b)>0$ we introduce the scaled point process $M^b$ (on $[0,\infty)$) by
$$
M^b(t):=N^b_+(\nu_+(b)^{-1}t),\quad t\ge 0.
$$
In each of Scenarios \ref{s1}--\ref{s3} described in
the Introduction  we will prove (under additional
technical assumptions) the convergence
\begin{align}\label{weak}
M^b \overset{d}{\longrightarrow} \Pi_\rho\quad \text{as $b\to\infty$},
\end{align}
where $\overset{d}{\longrightarrow}$ denotes weak convergence of point processes
(see e.g.\ \cite{K02}) under the probability measure $\BP$,
$\rho\in[0,1)$ is explicitly determined by the
characteristics of $(X_t)$ (see Theorem \ref{geocompound}),
and the geometrically compound Poisson process $\Pi_\rho$ was
defined in the Introduction. If $\rho=0$, then $\Pi_\rho$ is a unit rate
Poisson process. Actually we will prove the weak convergence
of $\BP_\sigma(M^b\in\cdot)$  for an essentially
arbitrary initial distribution $\sigma$.

In Scenarios 1 and 2 we assume that the jumps (from high enough levels)
in the respective processes are dominated in distribution.
This means that there exists a $u_0\in\R$ and a family of
non-increasing (right-continuous) functions
$(\overline{H}(u,\cdot))_{u\ge u_0}$ 
such that $\sup_{x\ge u}J(x, (z,\infty)) \le \overline{H}(u,z)$ for all $z\in\R$ and $u\ge u_0$.
In other words, denoting by $\xi(x)$ a generic r.v.\ with the distribution
$J(x,\cdot),$ this means that there exist r.v.'s $\oxi(u)$  
such that 
\begin{align}\label{Domino2}
\xi (x) \dle \oxi  (u), \qquad x\ge u \ge u_0.
\end{align}
We assume that $\BE \oxi(u)<\infty$. 
In Scenario 3 we will assume in addition that there exist r.v.'s $\uxi(u)$ such that
\begin{align}\label{Domino3}
\uxi  (u)\dle \xi (x), \qquad x\ge u \ge u_0.
\end{align}
Further, put
$$
\omu (u) := \sup_{x\ge u} \mu (x), \quad
\olambda (u) := \sup_{x\ge u} \lambda (x), \quad
\ulambda (u) := \inf_{x\ge u} \lambda (x).
$$

Next we will make Scenarios \ref{s1}--\ref{s3} more precise.

\begin{assumption}\label{as1}\rm
We have $\mu(y)\to -\infty$ as $y\to\infty$, and there exists a $u_0\in\R$ such that 
\eqref{Domino2} holds and $J(x,(-\infty,0])=0$ for all $x\ge u_0$. Moreover,
\begin{align}\label{NegDrift1}
\BE \oxi(u_0) + \omu(u_0)/\olambda(u_0)  <0.
\end{align}
\end{assumption}

\begin{assumption}\label{as2}\rm
We have $\lambda(y)\to  \infty$ as $y\to\infty$.
Furthermore, there is a $u_0\in\R$ such that \eqref{Domino2} holds,
$\mu(y)>0$ for all $y\ge u_0$,
$J(x,[(u_0-x)^+,\infty))=0$ for all $x\in\R$, and
\begin{align}\label{NegDrift2}
\BE\oxi(u_0) + \omu (u_0)/\ulambda(u_0)  <0.
\end{align}
\end{assumption}

\begin{assumption}\label{as3}\rm
As $y\to\infty$ we have $\mu(y)\to\mu (\infty)\in \R\setminus\{0\}$, 
$\lambda(y)\to \lambda(\infty)\in [0,\infty)$. 
There is some $u_0\in\R$ such that \eqref{Domino2} and \eqref{Domino3} hold,
$J(x,(-\infty,0))=0$ for all $x\ge u_0$ in case $\mu(\infty)<0$, and
$J(x,[(u_0-x)^+,\infty))=0$ for all $x\in\R$ in case $\mu(\infty)>0$.
Furthermore we have for $y\to\infty$ that $\oxi(y), \uxi(y) \dto \xi(\infty)$,
where $\xi(\infty)$ is an integrable  r.v.\ satisfying 
\begin{align}\label{NegDrift3}
\BE\xi( \infty) + \mu(\infty)/\lambda(\infty) <0.
\end{align}
\end{assumption}

\begin{remark}\label{r11}\rm
Each of the inequalities \eqref{NegDrift1}-\eqref{NegDrift3} 
implies the ergodicity condition \eqref{C7}. In case of
\eqref{NegDrift1} and \eqref{NegDrift2} this is due to the monotonicity
properties of $\oxi(u)$, $\omu$, $\olambda$, and $\ulambda$.
\end{remark}

To state the theorem we write $H$ for the set of all $x\in\R$ such
that $\BP_x(\tau(u)<\infty)=1$ for some $u\in \R$ satisfying
$\nu(u)>0$, where 
\begin{align}\label{tau(u)}
\tau(u):=\inf\{t>0:N^u(t)\ge 1\},\quad u\in\R,
\end{align} 
is the smallest point of $N^u$ and $\inf\emptyset:=\infty$.

\begin{theorem}\label{geocompound} Let one of Assumptions
\ref{as1}-\ref{as3} be satisfied and assume that $\nu(b)>0$ for all sufficiently large
$b$. Let $\sigma$ be a distribution on $\R$ that is supported by $H$.
Then $\BP_\sigma(M^b\in\cdot)$ converges weakly to $\BP(\Pi_\rho\in\cdot)$
as $b\to\infty$. The number $\rho$ is given by
$\rho=0$ in case of \ref{as1}, \ref{as2}, and in case \ref{as3} by
\begin{align}\label{Rho}
\rho =\begin{cases}
- \frac{\lambda (\infty)}{\mu (\infty)} \, \BE\xi(\infty), &\text{if $\mu(\infty)<0$},\\
 1- \frac{w \mu  (\infty)}{\lambda(\infty)}, &\text{if $\mu(\infty)>0$},
\end{cases}
\end{align}
where $w$ is the only positive number satisfying the equation
\begin{equation}\label{Eqn_w}
\BE e^{w\xi (\infty)} = 1- w\mu  (\infty)/\lambda (\infty).
\end{equation}
\end{theorem}

\begin{remark}\label{rem88}\rm
If $\lambda>0$ on $D_\mu$ then $\pi$ cannot be concentrated
on $D_\mu$, and \eqref{rice} implies the existence of a $u\in\R$ such that
$\nu(u)>0$. Lemma \ref{l3} below 
then implies that $\pi$-almost all $x\in\R$ belong to $H$.
\end{remark}

\begin{remark}\label{remf}\rm Let $g:\R\rightarrow\R$
be a strictly increasing continuously differentiable function  such that $g(x)\to\infty$
as $x\to\infty$. Then $(X^{g}_t):=(g(X_t))$ is again a piecewise-deterministic
Markov process as defined in Section \ref{sec2}. The characteristics
of $(X^g_t)$ are given by
$\mu^{g}(y) = g'(g^{-1}(y))\mu(g^{-1}(y))$, $\lambda^{g}(y)=\lambda(g^{-1}(y))$,
and $J^g(y,\cdot)=J(g^{-1}(y),g^{-1}(y+\cdot)-g^{-1}(y))$.
If the point processes of upcrossings defined in terms of $X^g$ satisfy
a compound limit theorem as in Theorem \ref{geocompound}, then so do the
corresponding processes defined in terms of $(X_t)$. 
Therefore the assertion of the theorem remains true in the more general case, when
one of the Assumptions \ref{as1}-\ref{as3} holds for the transformed
process $(X^g_t)$. 
\end{remark}

As a corollary we obtain that the first crossing time
$$
T(b) := \inf\{t>0:\, X_t \ge b\}
$$
is asymptotically exponentially distributed.

\begin{corollary}\label{c1} Under the assumptions of Theorem \ref{geocompound},
we have for any $s\ge 0$ that
\begin{align}\label{Expon_Law} 
\BP_\sigma((1-\rho) \nu(b) T (b)>s) \to e^{-s} \quad\text{as $b\to\infty$}.
\end{align}
\end{corollary}
{\sc Proof:} For any $s\ge 0$ and $b\ge u_0$ we have $\nu_+(b)=\nu(b)$  and
\begin{align}\label{67} \notag
\BP_\sigma((1-\rho)& \nu(b) T (b)>s)=\BP_\sigma(X_0<b,M^b((1-\rho)^{-1}s)=0)\\
&=\BP_\sigma(M^b((1-\rho)^{-1}s)=0)-\BP_\sigma(X_0\ge b,M^b((1-\rho)^{-1}s)=0).
\end{align}
The second term on the right-hand side of \eqref{67} converges to $0$ as $b\to\infty$.
As any fixed finite number of points (in our case $0$ and $(1-\rho)^{-1}s$)
are almost surely not contained in $\Pi_\rho$, we obtain from \eqref{weak}
and a standard property of weak convergence of point processes (see \cite{K02})
that the first term in \eqref{67} converges
to $\BP(\Pi_\rho((1-\rho)^{-1}s)=0)=e^{-s}$.\qed


\begin{remark}\label{r23}\rm
Define $T_1(b):=T(b)$ and, inductively,
$T_{n+1}(b):=\inf A_n$, $n\ge 1$, where $A_n$ is the set of all $t>T_n(b)$ such
that $X_t\ge b$ and $X_s<b$ for some $s\in(T_n(b),t)$. 
Under the assumptions of Theorem \ref{geocompound},
we obtain for any $n\ge 1$ and $s\ge 0$ as above that
\begin{align}\label{Gamma_Law} 
\BP_\sigma((1-\rho) \nu(b) T_n(b)>s) \to \BP(\Pi_\rho((1-\rho)^{-1}s)\le n-1)
 \quad\text{as $b\to\infty$}.
\end{align}
An easy calculation shows that, for instance,
\begin{align*}
\BP(\Pi_\rho((1-\rho)^{-1}s)\le 1)&=e^{-s}(1+(1-\rho)s),\\
\BP(\Pi_\rho((1-\rho)^{-1}s)\le 2)&=e^{-s}\Big(1+(1-\rho^2)s+\frac{(1-\rho)^3s^2}{2}\Big).
\end{align*}
\end{remark}

\begin{corollary}\label{c2} Let the assumptions of Theorem \ref{geocompound},
be satisfied and let $B\subset[0,\infty)$ be a bounded Borel set whose boundary
has Lebesgue measure $0$. Then
$M^b(B)\overset{d}{\longrightarrow}\zeta_B$ as $b\to\infty$, where
$\zeta_B$ is a non-negative integer-valued r.v.\ with the Laplace transform
\begin{align}\label{LT}
\BE\exp[-z\zeta_B]=\exp\Big[-|B|(1-\rho)\Big(1-\frac{1-\rho}{e^{z}-\rho}\Big)\Big],
\quad z\ge 0.
\end{align}
Here $|B|$ denotes the Lebesgue measure of $B$.
\end{corollary}
{\sc Proof:} The right-hand side of \eqref{LT} is just the Laplace transform
of $\Pi_\rho(B)$, see also the comments after \eqref{5.11}.
Hence the result is a direct consequence of Theorem \ref{geocompound} and
Theorem 16.16 in \cite{K02}.\qed

\begin{remark}\label{r65}\rm The random variable $\zeta_B\overset{d}{=}\Pi_\rho(B)$
is infinitely divisible with a L\'evy measure having the mass
$|B|(1-\rho)^2\rho^{k-1}$ at $k\ge 1$.
\end{remark}

Before proving Theorem \ref{geocompound} we will provide several lemmas. 
For $u\in\R$ we write $N^u(\infty):=\lim_{t\to\infty}N^u(t)$. 
We also recall definition \eqref{tau(u)}.

\begin{lemma}\label{l3} Assume that $u\in \R$ satisfies $\nu(u)>0$.
Then $\BP(\tau(u)<\infty)=1$ and $\BP(N^u(\infty)=\infty)=\BP_u(N^u(\infty)=\infty)=1$. 
Moreover, $\nu(u)=(\BE_u(\tau(u))^{-1}$.
\end{lemma}
{\sc Proof:} Take $x\in\R$.
The strong Markov property implies that $(X_t)_{t\le\tau(x)}$ and
$(\I\{\tau(x)<\infty\}X_{\tau(x)+t})_{t\ge 0}$ are independent
for any initial distribution. This fact will be often used in the sequel.
In particular, $N^x$ is a renewal process (with a possibly defective
distribution of interpoint distances).
Under $\BP$, $N^x$ is also a stationary point process. If
$\nu(x)>0$ this clearly implies that $\BP(\tau(x)<\infty)=1$.
The equation $\nu(x)=(\BE_x(\tau(x))^{-1}$ is then a consequence
of the elementary renewal theorem. In particular $\BE_x\tau(x)<\infty$, 
so that the equations $\BP(N^x(\infty)=\infty)=\BP_x(N^x(\infty)=\infty)=1$
are obvious.\qed 

\vspace{0.3cm}
For $u\in \R$ we define an increasing sequence $\tau_n(u)$, $n\ge 0$,
of stopping times inductively by $\tau_0(u):=0$ and
$\tau_{n+1}(u):=\inf\{t>\tau_n(u):N^u(t)\ge n+1\}$. Hence $\tau_1(u)=\tau(u)$
and $N^u(t)$ is the cardinality of $\{n\ge 1:\tau_n(u)\le t\}$.
If $\nu(u)>0$ then Lemma \ref{l3} implies for all $n\ge 1$ that
$\BP(\tau_n(u)<\infty)=\BP_u(\tau_n(u)<\infty)=1$.

\begin{lemma}\label{l2} Assume that $u\in \R$ satisfies $\nu(u)>0$
and let $b\in \R$. Then $\nu(b)>0$ iff $\BP_u(\tau(b)<\infty)>0$.
In this case $\BP_u(\tau(b)<\infty)=1$.
\end{lemma}
{\sc Proof:} Assume that $\nu(b)>0$. Since $\BP(\tau(u)<\infty)=1$ by
Lemma \ref{l3} we must have that $\BP_u(\tau(b)<\tau(u))>0$. Since
$\BP_u(N^u(\infty)=\infty)=1$ we can use a geometrical trial argument
to get $\BP_u(\tau(b)<\infty)=1$. 
Assume, conversely, that $\BP_u(\tau(b)<\infty)>0$. Then we must have
$\BP_u(\tau(b)<\tau(u))>0$ and hence $\BP_u(\tau(b)<\infty)=1$.
Lemma \ref{l3} implies that $\BP(\tau(b)<\infty)=1$. Therefore $N^b$ is 
a non-empty and stationary point process under $\BP$ and must hence
have a positive intensity $\nu(b)$.\qed

\vspace{0.3cm}
Our next lemma deals with the probabilities
\begin{align}
\gamma(u,b):=\BP_b(\tau(u)>\tau(b)),\quad u,b\in \R.
\end{align}

\begin{lemma}\label{l4} Let $u\in\R$ satisfy $\nu(u)>0$. Then we have for all
$b\in \R$ that
\begin{align}\label{4.5}
\BP_u(N^b(\tau(u))=k)&=\frac{\nu(b)}{\nu(u)}(1-\gamma(u,b))^2\gamma(u,b)^{k-1},\quad k\ge 1,\\
\label{4.6}
\BP_u(N^b(\tau(u))=0)&=1-\frac{\nu(b)}{\nu(u)}(1-\gamma(u,b)).
\end{align}
In particular,
\begin{align}\label{4.7}
\BE_uN^b(\tau(u))&=\frac{\nu(b)}{\nu(u)},\\
\label{4.72}
\BE_u\Big[1-\exp[-rN^b(\tau(u))]\Big]&=\frac{\nu(b)}{\nu(u)}
\Big(1-\gamma(u,b)-\frac{(1-\gamma(u,b))^2}{e^r-\gamma(u,b)}\Big),\quad r\ge 0.
\end{align}
\end{lemma}
{\sc Proof:} Equation \eqref{4.7}, i.e.\ $\nu(u)\BE_uN^b(\tau(u))=\nu(b)$
is an equilibrium equation that can be formulated for general stationary 
point processes. But in our case we can give a simpler argument as follows.
We assume that $\BP_u$ is the underlying probability measure. 
For any $n\ge 1$ the process $(X_{\tau_n(u)+t})_{t\ge 0}$ is Markov with
distribution $\BP_u((X_t)\in\cdot)$. This fact and the strong Markov property imply
that $X_n:=N^b(\tau_{n}(u))-N^b(\tau_{n-1}(u))$, $n\ge 1$, are i.i.d. 
>From the law of large numbers we obtain
$\BP_u$-a.s.\ that $n^{-1}\sum^n_{k=1}X_n\to \BE_uX_1=\BE_uN^b(\tau(u))$ as
$n\to\infty$. Assume that $\nu(b)>0$. From Lemma \ref{l2} we have
$\BP_u(\tau(b)<\infty)=1$ so that we can use the
laws of large numbers for independent random variables and renewal processes
to obtain $\BP_u$-a.s.\ that
$$
\frac{1}{n}\sum^n_{k=1}X_n=\frac{\tau_n(u)}{n}\frac{1}{\tau_n(u)}N^b(\tau_n(u))
\to \frac{\BE_u\tau(u)}{\BE_b\tau(b)}=\frac{\nu(b)}{\nu(u)}\quad \text{as $n\to\infty$},
$$
where we have again used Lemma \ref{l3}. This implies \eqref{4.7}.
In case $\nu(b)=0$ we have $\BP_u(\tau(b)<\tau(u))=0$, so that \eqref{4.7} is valid as well.

Next we use the strong Markov property to obtain for $k\ge 1$
\begin{align*}
\BP_u(N^b(\tau(u))=k)&=\BP_u(\tau(b)<\tau(u),N^b(\tau(u))=k)\\
&=\BE_u\I\{\tau(b)<\tau(u)\}\BP_b(N^b(\tau(u))=k-1)
=p\gamma(u,b)^{k-1}(1-\gamma(u,b)),
\end{align*}
where $p:=\BP_u(\tau(b)<\tau(u))=\BP_u(N^b(\tau(u))>0)$. If $\gamma(u,b)<1$, then we have,
in particular, that
$$
\BE_uN^b(\tau(u))=\frac{p}{1-\gamma(u,b)}.
$$
Comparing this with \eqref{4.7} yields $p=(1-\gamma(u,b))\nu(b)/\nu(u)$ and hence
\eqref{4.5} and \eqref{4.6}. In case $\gamma(u,b)=1$ these relations are true as well. 
Equation \eqref{4.72} follows from a direct computation.
\qed

\vspace*{0.3cm}
The assumptions of Theorem \ref{geocompound} are used in the following 
key lemma.

\begin{lemma}\label{l6}
If one of the Assumptions \ref{as1} or \ref{as2} is met, then
\begin{equation}\label{Goes_Down1}
\lim_{b\to\infty}\gamma(u,b)= 0
\end{equation}
for all sufficiently large $u$. If Assumption \ref{as3} is satisfied, then
\begin{equation}\label{Goes_Down2}
\lim_{u\to\infty}\liminf_{b\to\infty}\gamma(u,b)=\lim_{u\to\infty}\limsup_{b\to\infty}\gamma(u,b)
=\rho \in(0,1),
\end{equation}
where $\rho$ is defined in \eqref{Rho}.
\end{lemma}
{\sc Proof:} Without loss of generality, we can assume that $\overline{H}(\cdot,z)$ is 
non-increasing for any
$z\in\R$ (or, equivalently, that $\oxi(u)\dle \oxi(v)$ for $u > v$)
and that $\underline{H}(\cdot,z)$ is non-decreasing.
In the whole proof we will assume that $u_0$ is chosen according to
one of the Assumptions \ref{as1}-\ref{as3}. 
Note that the drift condition
\eqref{NegDrift1} (resp.\ \eqref{NegDrift2}) holds with $u$ in place of $u_0$.
It is then no loss of generality to assume that $\mu(u)\ne 0$ for all $u\ge u_0$.
We always take $u,b\in\R$ such that
$b> u\ge u_0$. In both cases \ref{as1} or \ref{as2} the argument will run roughly as follows. Due 
to the imposed conditions, for a large enough initial value $b$, the trajectory of 
$(X_t)$ will very quickly drop by a given large quantity $C$. Since in the part 
of the state space above the level $u$ the process can be shown 
to be dominated (at its jump points) by a random
walk with i.i.d.\ jumps and a negative trend, we can choose $C$ large enough to ensure
that the process will not climb back by $C$ prior to dropping below the level $u$.
If $(X_t)$ does not drop quickly enough from a high level, then it is likely there will
be several crossings of that level before the process returns to the range of its
`normal values'. This case requires the more restrictive conditions 
formulated in Assumption \ref{as3}.

First assume that Assumption \ref{as1} is met. Fix an arbitrary $\varepsilon >0$. 
Since the process
has a negative drift in the half-line $[u,\infty)$, it can only exceed the level $b>u$ by
a jump, so we can restrict ourselves to considering the values $X_{T_1},
X_{T_2},\ldots$: 
\begin{equation}\label{U_bound1}
\BP_x \bigl(\tau(b)<\tau(u)\bigr) \le \BP_x\bigl(\sup\{X_{T_k} :\, T_k \le \tau(u)\} \ge b\bigr),
\quad x\ge u.
\end{equation}
Further, for $x,t$ such that $q(x,t)>u$ we have
\begin{align}\label{4.81}
\BP_x(T_1 >t) =\exp \bigg[- \int_0^t \lambda (q(x,s)) ds\bigg]
\ge\exp[-\olambda (u)t] =\BP \bigl(\tau (\olambda (u))>t\bigr),
\end{align}
where  $\tau(w)$ is a  r.v.\   following the exponential  distribution with
parameter~$w$. Therefore one can easily see that the right-hand side  of
\eqref{U_bound1} does not exceed $ \BP (S \ge b-x),$ where $S:=\sup_{k\ge 1} S_k$ is the
global supremum of a random walk
\begin{equation}\label{RW}
S_k=\zeta_1+\cdots + \zeta_k, \qquad k\ge 1,
\end{equation}
with i.i.d.\ jumps $\zeta_k\deq\oxi(u)+\omu(u)\tau(\olambda(u))$, where $\oxi(u)$
and $\tau (\olambda (u))$ are independent of each other. Since $\BE\zeta_k<0$
by \eqref{NegDrift1}, $S$ is a proper r.v., and we can choose $C$ so large 
that $\BP(S\ge C)<\varepsilon$.

Next we assume $b\ge u+C$. For any $t\ge 0$ the equation $q(b,t)=b-C$ has a unique
solution $t=t(b,C)$. Since $\mu(y)\to-\infty$ as $y\to\infty$, we have 
$t(b,C)\to 0$ as $b\to\infty$. In particular, we obtain from \eqref{4.81} that
$\BP_b(T_1 \le t(b,C))<\varepsilon$ for all large enough $b$.
Then we have $\BP_b (X_{t(b,C)}=b-C)>1-\varepsilon$,
and finally, due to \eqref{U_bound1} and our choice of $C$, that 
$\BP_b(\tau(b)<\tau(u))<2\varepsilon$. Since $\varepsilon$ was arbitrary small, 
this completes the proof of the lemma in the first case.

Now suppose that Assumption \ref{as2} holds. In this case, jumps from levels 
$x\ge u$ are negative, and
we can concentrate on the values $X_{T_1-}, X_{T_2-},\ldots$. For a given 
$\varepsilon >0$ choose an $C<\infty$ such that 
$\BP(S \ge C)<\varepsilon$ for the random walk \eqref{RW} with
i.i.d.\ jumps $\zeta_k\deq \oxi (u) +\omu (u) \tau(\ulambda(u))$, where $\oxi(u)$ and
$\tau (\olambda (u))$ are independent of each other. 
Consider a stopping time $T$ with values in $\{T_1,T_2,\ldots\}$.
Then, as one can easily see, given
that $X_{T-}<b-C$, the probability of the process exceeding $b$ 
on the time interval $[T,\infty)$ prior to dropping
below $u$ will again be less than $\varepsilon$.

Since the deterministic drift is now positive on $(u_0,\infty)$, 
$$
\BP_x ( T_1 >t) =\exp \biggl\{ - \int_0^t \lambda (q(x,s)) ds\biggr\}
\le \exp\{- \ulambda (x)t\} =\BP \bigl(\tau (\ulambda(x))>t\bigr),\quad x\ge u.
$$
Therefore, given $X_0=b$, one has $X_{T_1 -}\dle b+\omu (u)\tau (\ulambda (b))$, and it is
not difficult to see that, for $m\ge 1$, $\overline X_m := \sup_{T_1\le t \le T_m} X_t,$
$\underline X_m := \inf_{T_1\le t \le T_m} X_t,$
\begin{multline}
\BP_b \bigl( \{\overline X_m \ge  b\}
 \cup \{ \underline X_m \ge b-C\} \bigr)\\
\le \BP \Bigl(  \max_{1\le k\le m} S_k \ge -\omu (u)\tau (\ulambda (b-C))\Bigr)
+\BP (S_m \ge -C),
\label{OverX}
\end{multline}
where $(S_k)$ is a random walk given by \eqref{RW} with i.i.d.\ jumps $\zeta_k\deq \oxi
(u) +\omu (u) \tau (\ulambda (b -C))$, and $(S_k)$, $\tau (\ulambda (b-C))$ that
appear together  under the probability sign in \eqref{OverX} are independent of each
other. Now choose $m$ so large that $\BP (S_m \ge -C)<\varepsilon$ (this is possible due
to \eqref{NegDrift2}), and then $b$ so large that the first term on the right-hand side
of \eqref{OverX} is also less than $\varepsilon$. The latter is possible due to the following
observation. Setting
$$
S'_k:=\xi'_1+\cdots + \xi'_k, \quad
S''_k:=\xi''_1+\cdots + \xi''_k, \quad k\ge 1,
$$
where $(\xi'_k)$ and $(\xi''_k)$ are independent sequences of i.i.d.\ r.v.'s with
$\xi'_k\deq \oxi (u)$, $\xi''_k\deq \tau (\ulambda(b-C))$, the event in that term 
is contained in
$$
\max_{1\le k\le m} S'_k \ge -\omu(u) \max_{1\le k\le m} S''_k -\omu(u)\xi''_{m+1}.
$$
The  r.v.\ on the left-hand side is a.s.\ negative, with a
distribution independent of $b$. Because it is assumed that $\lambda(y)\to\infty$
as $y\to\infty$, the distribution of the right-hand side
converges to $\delta_0$ as $b\to\infty$.

Thus, on the event complementary to the one on the left-hand side of \eqref{OverX}, the
process $(X_t)$ will drop at one of the times $T_1,\ldots, T_m$ below the level $b-C$
(denote this epoch by $T^*$), without having continuously crossed the level $b$ prior to
that time. Also, due to our choice of~$C$ and to the strong Markov property, the process
will reach the level~$b$ on the time interval $[T^*,\tau(u)]$ with probability less
than~$\varepsilon$. This means that  $\BP_b(\tau(b)<\tau(u))< 3\varepsilon$ and hence proves 
the lemma in the case when Assumption \ref{as2} holds.

Now consider the case when Assumption \ref{as3} holds.
Assume first that $\mu(\infty)<0.$ Then crossing the level
$b$ can only occur due to a jump, and since to get from a level $x>b$ down to level $u$
will require a continuous downcrossing of $b$, we obtain that
\begin{equation}
\BP_b(\tau(b)<\tau(u)) = \BP_b\bigl(\sup\{ X_{T_k}- b:\, T_k<\tau(u) \}>0\bigr).
\label{UV}
\end{equation}
Next we observe that, for the segment of the process in the time interval $[0,\tau(u)]$, one
has $\uS_k \dle X_{T_k}-b \dle \oS_k,$ where $(\oS_k)_{k\ge 1}$ and $(\underline
S_k)_{k\ge 1}$ are random walks with i.i.d.\ jumps
$$
\oxi_k\deq \oxi (u) + \omu (u) \tau (\olambda(u)),\qquad
\uxi_k\deq \uxi (u) + \umu (u) \tau (\ulambda(u)),
$$
respectively, where we again make the usual independence assumptions.
Due to \eqref{NegDrift3} and uniform integrability of $(\oxi (u))$, we get
$\BE \uxi_k\le \BE \oxi_k<0$ for all large enough $u$, so that then
$$
\uS^u:= \sup_{k\ge 1} \uS_k \dle \oS^u:= \sup_{k\ge 1} \oS_k<\infty\quad\text{ a.s.}
$$
It is not difficult to see that, for $b>2u$,
$$
\BP(\uS^u >0) +R(u,b)
\le \BP_b\bigl(\sup\{ X_{T_k}- b:\, T_k<\tau(u)\}>0\bigr)
\le \BP (\oS^u >0),
$$
where
$$
R(u,b):=\BP(\sup\{\uS_k:k\le\eta\}>0)-\BP(\uS^u>0)
-\BP(\umu (u)\tau (\ulambda (u))>b/2-u),
$$
and $\eta:=\inf\{k>0:\uS_k<-b/2\}$. Since clearly $\eta\to\infty$ a.s.\  as $b\to\infty$,
we obtain that
\begin{align}\label{claim}
\lim_{b\to\infty}R(u,b)=0.
\end{align}

By virtue of Theorem 6 and Condition $B$ on p.114 of \cite{Bo76}, we obtain that  
$$
\lim_{u\to\infty}\BP(\uS^u>0)=\lim_{u\to\infty}\BP(\oS^u >0)=\BP(S>0),
$$ 
where $S=\sup_{k\ge 1} S_k$ for a random walk
with i.i.d.\ jumps $\zeta_k\deq \xi (\infty)+\mu(\infty)\tau(\lambda(\infty))$. 
Because $(S_k)$ has a negative drift, it is well-known that
$\BP(S >0)$ is given as in \eqref{Rho}, see e.g.\ Theorem VIII.5.7 and Corollary III.6.5
in \cite{As03}.

The argument in the case when $\mu (\infty) >0$ is very similar, with the main
difference being the value of $\BP(S>0)$. But again it is well-known that
this value is given as in \eqref{Rho}, see e.g.\ Theorem X.5.1 in \cite{As03}.\qed

\begin{lemma}\label{l7} Assume that
one of the Assumptions \ref{as1}-\ref{as3} is satisfied and that
$\nu(u)>0$ for all sufficiently large $u$. Then we have,
for any $\delta>0$, 
$$
\lim_{u\to\infty}\limsup_{b\to\infty}
\frac{1}{\nu(b)}\BE_u\int^\infty_0\I\{\tau(u)\ge s>\delta/\nu(b)\}N^b_+(ds)=0.
$$
\end{lemma}
{\sc Proof:} We take $u_0\in\R$ according to Assumptions \ref{as1}-\ref{as3}.
and assume that $\mu(u)\ne 0$ and $\nu(u)>0$ for all $u\ge u_0$. The numbers $u,b$ are always
chosen so that $b>u\ge u_0$. We first note that
\begin{align}\label{4.51}
\BP_u(N^b_+(\tau(u))=N^b(\tau(u)))=1,\quad u\ge u_0.
\end{align} 
In case $\mu<0$ on $(u_0,\infty)$ this is due to the absence of negative jumps from
a level above $b$. In case $\mu>0$ on $(u_0,\infty)$ we even have $N^b_+=N^b$
because there is no positive jump to a level above $b$.
Next we define the stopping times $\tau^+_k(b)$, $k\ge 1$, in terms of
of $N^b_+$ as $\tau_k(b)$ in terms of $N^b$.
Then we have for any integer $m\ge 1$ that
$$
\frac{1}{\nu(b)}\int^\infty_0\I\{\tau(u)\ge s>\delta/\nu(b)\}N^b_+(ds)
=\eta_m(u,b)+\zeta_m(u,b),
$$
where 
\begin{align*}
\eta_m(u,b)&:=\frac{1}{\nu(b)}\sum^m_{k=1}\I\{\tau(u)\ge \tau^+_k(b)>\delta/\nu(b)\},\\
\zeta_m(u,b)&:=\frac{1}{\nu(b)}\int^\infty_0\I\{\tau(u)\ge s>\delta/\nu(b), s>\tau^+_m(b)\}N^b_+(ds),
\end{align*}
Under $\BP_u$ we have $\{\tau(u)\ge\tau^+_1(b)\}\subset \{\tau(u)>\tau(b)\}$.
This follows as at \eqref{4.51}. Hence we get
\begin{align*}
\BE_u\eta_m(u,b)&\le \frac{m}{\nu(b)}\BE_u\I\{\tau(u)>\tau(b)>\delta/\nu(b)\}\\
&\le \frac{m}{\delta}\BE_u\I\{\tau(u)>\tau(b)\}\tau(u)
\end{align*}
by Markov's inequality. Furthermore,
$$
\BP_u(\tau(u)>\tau(b))=\BP_u(N^b(\tau(u))>0)\le\BE_u N^b(\tau(u))=\frac{\nu(b)}{\nu(u)},
$$
where the last equality comes from \eqref{4.7}. 
Because $\BE_u\tau(u)=1/\nu(u)<\infty$
and $\nu(b)\to 0$ as $b\to\infty$, we can use dominated convergence to conclude,
for any fixed $m\ge 1$, that $\BE_u\eta_m(u,b)\to 0$ as $b\to\infty$.

To deal with $\zeta_m(u,b)$ we use the simple estimate
\begin{align*}
\BE_u\zeta_m(u,b)&\le \frac{1}{\nu(b)}\BE_u\I\{N^b_+(\tau(u))\ge m+1\}N^b_+(\tau(u))\\
&= \frac{1}{\nu(b)}\BE_u\I\{N^b(\tau(u))\ge m+1\}N^b(\tau(u))
\end{align*}
and \eqref{4.5} to obtain
\begin{align*}
\BE_u\zeta_m(u,b)\le \frac{(1-\gamma(u,b))^2}{\nu(u)}\sum^\infty_{k=m+1}k\gamma(u,b)^{k-1}
=\frac{\gamma(u,b)^m}{\nu(u)}(m(1-\gamma(u,b))+1),
\end{align*}
where the equality comes from a direct calculation.
Let $\varepsilon>0$. By Lemma \ref{l6} we then find an $m\ge 1$ such
that $\limsup_{b\to\infty}\BE_u\zeta_m(u,b)\le \varepsilon$ as soon as $u$ is sufficiently
large. Together with the first part of the proof this implies the assertion.
\qed

\vspace*{0.3cm}
{\sc Proof of Theorem \ref{geocompound}:} We first prove the result in the
stationary case, i.e.\ we take $\sigma:=\pi$.
Let $f:[0,\infty)\rightarrow[0,\infty)$ be a continuous function with compact support
and 
$$
L(f,b):=\BE\exp\Big[-\int f(\nu(b)s)N^b_+(ds)\Big],\quad b\in\R.
$$
We will show that
\begin{align}\label{5.11}
\lim_{b\to\infty} L(f,b)
=\exp\Bigg[-(1-\rho)\int^\infty_0 \bigg(1-\frac{1-\rho}{e^{f(t)}-\rho}\bigg)dt\Bigg].
\end{align}
The right-hand side of \eqref{5.11} coincides with
$\BE\exp\big[-\int f(s)\Pi_\rho(ds)\big]$, as can easily be confirmed by using
Lemma 12.2 (i),(iii) in \cite{K02}. Hence Theorem 16.16 in
\cite{K02} implies the assertion \eqref{weak}.

We assume now that $u_0\in\R$ has been chosen according to one of the
Assumptions \ref{as1}-\ref{as3}.
Without loss of generality we can also assume that $u\notin D_\mu$ and 
$\nu(u)>0$ for all $u\ge u_0$. 
In the following we will always pick $b>u\ge u_0$.
We use the notation introduced before Lemma \ref{l4}. 
By the strong Markov property, the restrictions of $N^b_+$ to the 
(random) intervals $(\tau_n(u),\tau_{n+1}(u)]$, $n\ge 0$, are independent. Therefore,
\begin{align}\notag
L(f,b)&=\prod^\infty_{n=0}\BE\exp\Big[-\int_{(\tau_n(u),\tau_{n+1}(u)]} f(\nu(b)s)N^b_+(ds)\Big]\\
\label{4.76}
&=L_0(u,b)\prod^\infty_{n=1}\int\BE_u\exp\Big[-\int_{(0,\tau(u)]} f(\nu(b)(s+t))N^b_+(ds)\Big]
\BP(\tau_n(u)\in dt),
\end{align}
where the second equality is again a consequence of the strong
Markov property, and
$$
L_0(u,b):=\BE\exp\Big[-\int_{(0,\tau(u)]} f(\nu(b)s)N^b_+(ds)\Big].
$$
We claim that $\BP(N^b_+(\tau(u))>0)\to 0$
as $b\to\infty$, so that dominated convergence implies that
\begin{align}\label{R0}
R_0(u,b):=-\ln L_0(u,b)\to 0\quad\text{as $b\to\infty$}.
\end{align}
To prove the claim, we pick numbers $\varepsilon>0$ and 
$s>0$ to obtain that
\begin{align*}
\BP(N^b_+(\tau(u))>0)&=\BP(\tau^+_1(b)\le \tau(u))\\
&=\BP(\tau^+_1(b)\le \tau(u),\tau(u)>s)+\BP(\tau^+_1(b)\le \tau(u),\tau(u)\le s)\\
&\le \BP(\tau(u)>s)+\BP(\tau^+_1(b)\le s).
\end{align*}
For large enough $s$ the first term is smaller than $\varepsilon$. 
For the second term we have
$$
\BP(\tau^+_1(b)\le s)=\BP(N^b_+(s)>0)\le \BE N^b_+(s)=\nu(b)s.
$$
The right-hand side is getting smaller than $\varepsilon$ for all large enough $b$.

Defining
$$
h(u,b,t):=\BE_u\exp\Big[-\int_{(0,\tau(u)]} f(\nu(b)(s+t))N^b_+(ds)\Big]
$$
we obtain from \eqref{4.76} that
\begin{align*}
-\log L(f,b)&=-\sum^\infty_{n=1}\log \BE h(u,b,\tau_n(u))+R_0(u,b)\\
&=\sum^\infty_{n=1}(1-\BE h(u,b,\tau_n(u)))+R_0(u,b)+\sum^\infty_{n=1}\theta(1-\BE h(u,b,\tau_n(u))),
\end{align*}
where $|\theta(r)|\le cr^2$ for some (universal) constant $c>0$.
Using Campbell's theorem for the stationary point process $N^u$ 
(see e.g.\ equation (1.2.18) in \cite{BB}) gives
\begin{align}\label{log}
-\log L(f,b)=\nu(u)\int^\infty_0(1-h(u,b,t))dt+R_0(u,b)+R_1(u,b),
\end{align}
where the remainder term $R_1$ is defined by
\begin{align}\label{5.17}
R_1(u,b):=\sum^\infty_{n=1}\theta(1-\BE h(u,b,\tau_n(u))).
\end{align}

Using Jensen's inequality, the inequality $(1-e^{-x})\le x$, $x\ge 0$, 
and Campbell's theorem again we get
\begin{align*}
|R_1(u,b)|&\le c\sum^\infty_{n=1}(1-\BE h(u,b,\tau_n(u)))^2
\le c\sum^\infty_{n=1}\BE(1- h(u,b,\tau_n(u))^2\\
&\le c\sum^\infty_{n=1}\int^\infty_0 \Big(\BE_u\int_{(0,\tau(u)]} 
f(\nu(b)(s+t))N^b_+(ds)\Big)^2 \BP(\tau_n(u)\in dt)\\
&=c\nu(u)\int^\infty_0 \Big(\BE_u\int_{(0,\tau(u)]} f(\nu(b)(s+t))N^b_+(ds)\Big)^2 dt\\
&=c\frac{\nu(u)}{\nu(b)}\int^\infty_0 \Big(\BE_u\int_{(0,\tau(u)]} f(\nu(b)s+r)N^b_+(ds)\Big)^2 dr\\
&=c\frac{\nu(u)}{\nu(b)}\int^\infty_0 \Big(\int_0^\infty f(\nu(b)s+r)m_{u,b}(ds)\Big)^2 dr,
\end{align*}
where the measure $m_{u,b}$ is given by
$$
m_{u,b}(\cdot):=\BE_u\int^\infty_0\I\{s\in\cdot,0<s\le\tau(u)\}N^b_+(ds).
$$
By \eqref{4.7} the measure $m^*_{u,b}:=\frac{\nu(u)}{\nu(b)}m_{u,b}$ has total
mass $1$. Hence we can use Jensen's inequality to obtain that
\begin{align*}
|R_1(u,b)|&\le c\frac{\nu(b)}{\nu(u)}\int^\infty_0 \Big(\int_0^\infty f(\nu(b)s+t)m^*_{u,b}(ds)\Big)^2 dt\\
&\le c\frac{\nu(b)}{\nu(u)}\int^\infty_0\int_0^\infty f(\nu(b)s+t)^2 m^*_{u,b}(ds) dt.
\end{align*}
By Fubini's theorem and a change of variables
\begin{align}\label{rem}
|R_1(u,b)|\le c\frac{\nu(b)}{\nu(u)}\int^\infty_0\int_0^\infty \I\{t\ge \nu(b)s\}f(t)^2 dt\, m^*_{u,b}(ds)
\le c\frac{\nu(b)}{\nu(u)}\int^\infty_0 f^2(t)dt.
\end{align}

The main term in \eqref{log} equals
\begin{align}\notag
\nu(u)\int^\infty_0&(1-h(u,b,t))dt
=\frac{\nu(u)}{\nu(b)}
\int^\infty_0\BE_u\bigg[1-\exp\Big[-\int_{(0,\tau(u)]} f(\nu(b)s+t)N^b_+(ds)\Big]\bigg]dt\\
\notag
=&\frac{\nu(u)}{\nu(b)}
\int^\infty_0\BE_u\Big[1-\exp\big[-f(t)N^b_+(\tau(u))\big]\Big]dt+R_2(u,b)\\
\label{5.23}
=&(1-\gamma(u,b))\int^\infty_0 \bigg(1-\frac{1-\gamma(u,b)}{e^{f(t)}-\gamma(u,b)}\bigg)dt
+ R_2(u,b).
\end{align}
where (recall \eqref{4.51})
$$
R_2(u,b):=\frac{\nu(u)}{\nu(b)}
\int^\infty_0\BE_u\Big[\exp\big[-f(t)N^b(\tau(u))\big]-
\exp\Big[-\int_{(0,\tau(u)]} f(\nu(b)s+t)N^b_+(ds)\Big]\Big]dt,
$$
and we have used \eqref{4.72} to obtain the last equality.

To deal with the remainder term $R_2(u,b)$, we use the inequality
$$
\Big|\prod^n_{i=1} z_i-\prod^n_{i=1} w_i\Big|\le\sum^n_{i=1}|z_i-w_i|,
$$
for numbers $z_1,\ldots,z_n,w_1,\ldots,w_n$ of absolute value less than or equal to $1$.
This yields in case $X_0=u$ 
\begin{align*}
\Bigg|\exp\big[&-f(t)N^b(\tau(u))\big]-
\exp\Big[-\int_{(0,\tau(u)]} f(\nu(b)s+t)N^b_+(ds)\Big]\Big]\Bigg|\\
&\le \int_{(0,\tau(u)]}\big|\exp[-f(t)]- \exp[-f(\nu(b)s+t)]\big|N^b_+(ds).
\end{align*}
Hence we obtain for any $\delta>0$ that
\begin{align}\notag
|&R_2(u,b)|\le \frac{\nu(u)}{\nu(b)}
\BE_u\int^\infty_0\int_{(0,\tau(u)]}\I\{\nu(b)s\le\delta\}\big|\exp[-f(t)]
- \exp[-f(\nu(b)s+t)]\big|N^b_+(ds)dt\\
\label{5.33}
&+\frac{\nu(u)}{\nu(b)}\BE_u\int^\infty_0\int_{(0,\tau(u)]}\I\{\nu(b)s>\delta\}
\big|\exp[-f(t)]- \exp[-f(\nu(b)s+t)]\big|N^b_+(ds)dt.
\end{align}
As $b\to\infty$ we can use the uniform continuity of $f$ and \eqref{4.7}
to make the first term  arbitrarily small just by choosing $\delta$ small enough.
The second term \eqref{5.33} is smaller than
\begin{align*}
\frac{\nu(u)}{\nu(b)}&\BE_u\int^\infty_0\int_{(0,\tau(u)]}\I\{\nu(b)s>\delta\}
\big|1- \exp[-f(\nu(b)s+t)]\big|N^b_+(ds)dt\\
&+\frac{\nu(u)}{\nu(b)}\BE_u\int^\infty_0\int_{(0,\tau(u)]}\I\{\nu(b)s>\delta\}
\big|1-\exp[-f(t)]\big|N^b_+(ds)dt\\
\le&
\frac{\nu(u)}{\nu(b)}\BE_u\int^\infty_0\int_{(0,\tau(u)]}\I\{\nu(b)s>\delta\}
f(\nu(b)s+t)N^b_+(ds)dt\\
&+\frac{\nu(u)}{\nu(b)}\BE_u\int^\infty_0\int_{(0,\tau(u)]}\I\{\nu(b)s>\delta\}
f(t)N^b_+(ds)dt\\
\le&\frac{2\nu(u)}{\nu(b)}\Big(\int^\infty_0 f(t)dt\Big) \BE_u\int_{(0,\tau(u)]}\I\{\nu(b)s>\delta\}
N^b_+(ds).
\end{align*}
Hence we conclude from Lemma \ref{l7} that
\begin{align}\label{aha}
\lim_{u\to\infty}\limsup_{b\to\infty}|R_2(u,b)|=0.
\end{align}

Summarizing \eqref{log} and \eqref{5.23} gives
$$
-\log L(f,b)=(1-\gamma(u,b))\int^\infty_0 \bigg(1-\frac{1-\gamma(u,b)}{e^{f(t)}-\gamma(u,b)}\bigg)dt
+R_0(u,b)+R_1(u,b)+R_2(u,b).
$$
>From Lemma \ref{l6},  \eqref{R0}, \eqref{rem}, and \eqref{aha} we obtain \eqref{5.11}
and hence the assertion of the theorem in case $\sigma=\pi$.

Coupling is a well-known and elegant method to extend limit theorems beyond the
stationary setting. As we have not assumed ergodicity it is not possible to use
exact coupling as in Theorem 10.27 (i) in \cite{K02}. And the
shift-coupling assertion (ii) of
that theorem does not seem to be sufficient for our goals. So our strategy is to
use Thorisson's shift-coupling of point processes, see Lemma 11.7 in \cite{K02}.
Unless started otherwise we are working under the stationary probability
measure $\BP$. In a first step we extend
$(X_t)_{t\ge 0}$ to a stationary process $X:=(X_t)_{t\in\R}$,
such that the extended process is still right-continuous with
left-hand limits. We refer here to \cite{BB} for more details.
For $u\in\R$ we introduce as before the point process $N^u_+$ on $\R$. 
As usual we are identifying
a point process on $\R$ with a random (counting) measure on $\R$.
The scaled point process $M^u$ is defined by
$M^u(B):=N^u_+(\nu_+(u)^{-1}B)$ for any Borel set $B\subset\R$.
Also $\Pi_\rho$ can be extended to a stationary point process 
$\Pi'_\rho$ on $\R$. 
By stationarity it is then immediate that the weak convergence \eqref{weak} 
extends to $\R$. 

Next we introduce the space $\mathbf{D}$ of all mappings $z=(z_t)_{t\in\R}:\R\rightarrow\R$
that are right-continuous with left-hand limits equipped with
the $\sigma$-field $\mathcal{D}$ generated by the Skorohod topology
(see e.g.\ Theorem A2.2 in \cite{K02}). For any $s\in\R$ we define
the shift $\theta_s:\mathbf{D}\rightarrow\mathbf{D}$ by
$\theta_sz:=(z_{t+s})_{t\in\R}$. The distribution $\BP':=\BP(X\in\cdot)$ is
stationary, i.e.\ invariant under all these shifts. Let
$\mathcal{I}\subset \mathcal{D}$ denote the invariant $\sigma$-field, i.e.\
the system of all sets $A\in\mathcal{D}$ satisfying $\theta_sA=A$ for all
$s\in\R$. From the ergodic theorem (see Corollary 10.9 and Exercise 10.6 in \cite{K02})
we have for all
bounded and measurable $f:\mathbf{D}\rightarrow\R$ that
$(2t)^{-1}\int^t_{-t}f(\theta_sX)ds$ converges $\BP$-almost surely
to $\BE[f(X)|X^{-1}\mathcal{I}]$ as $t\to\infty$. 
The regenerative structure of $X$ implies
on the other hand that this limit must be a.s.\ constant.
Hence $f(X)$ and $X^{-1}\mathcal{I}$ are independent. In particular,
the $\sigma$-field $X^{-1}\mathcal{I}$ is a.s.\ trivial, so that $\BP'$
is ergodic (in the sense of ergodic theory).

Let $u\in\R$ have $\nu(u)>0$ and introduce the probability
measure
\begin{align}\label{Palm}
\BP^0_u(A):=\nu(u)^{-1}\BE\int^1_0\I\{\theta_sX\in A\}N^u(ds),\quad A\in\mathcal{D}.
\end{align}
(This is nothing but the {\em Palm probability measure} of $N^u$.)
A conditioning w.r.t.\  $X^{-1}\mathcal{I}$ shows that $\BP^0_u$ is also
trivial on $\mathcal{I}$.
By Lemma 11.7 in \cite{K02} we can hence assume without loss of generality that
there is a real-valued random variable $\tau$ satisfying
$\BP(\theta_\tau X\in\cdot)=\BP^0_u$. For any $b\in\R$ the scaled
point process of upcrossings of the level $b$ can be written as a measurable function
$M^b\equiv M^b(X)$ of $X$. Using stationarity it is easy to derive the weak convergence
of $\BP(M^b(\theta_\tau X)\in\cdot)$ to $\BP(\Pi'_\rho\in\cdot)$
from the weak convergence proved above. Let $X^+(z):=(z_t)_{t\ge 0}$ be
the restriction of the function $z$ on $[0,\infty)$. Since
$\BP(0\in\Pi'_\rho)=0$ we have weak convergence of
$\BP^0_u(M^b(X^+)\in\cdot)$ to $\BP(\Pi_\rho\in\cdot)$. From the strong Markov property
and \eqref{Palm} we have on the other hand that
$\BP^0_u(X^+\in\cdot)=\BP_u((X_t)_{t\ge 0}\in\cdot)$.
Hence we conclude the assertion for $\sigma=\delta_u$.

Finally we will prove the assertion for $\sigma=\delta_x$, where
$x\in\R$ satisfies $\BP_x(\tau(u)<\infty)=1$ for some $u\in \R$ with
$\nu(u)>0$. This is enough to conclude the theorem.
By the last assertion of Lemma \ref{l4} we can assume that $u>u_0$.
For $f$ as in \eqref{5.11} we have
\begin{align*}
\BE_x\exp\Big[-\int^\infty_0 f(\nu(b)s)N^b_+(ds)\Big]=L(b)-R(b),
\end{align*}
where 
\begin{align*}
L(b)&:=\BE_x\exp\Big[-\int^\infty_{\tau(u)} f(\nu(b)s)N^b_+(ds)\Big],\\
R(b)&:=\BE_x\Big(1-\exp\Big[-\int_0^{\tau(u)} f(\nu(b)s)N^b_+(ds)\Big]\Big)
\exp\Big[-\int^\infty_{\tau(u)} f(\nu(b)s)N^b_+(ds)\Big].
\end{align*}
By the strong Markov property,
\begin{align*}
L(b)&=\int\BE_u\exp\Big[-\int^\infty_{0} f(\nu(b)(s+t))N^b_+(ds)\Big]\BP_x(\tau(u)\in dt)\\
&=\int\BE_u\exp\Big[-\int^\infty_{0} f(s+\nu(b)t)M^b(ds)\Big]\BP_x(\tau(u)\in dt).
\end{align*}
Since $\nu(b)\to 0$ as $b\to\infty$ we can use the continuous mapping theorem
(see Theorem 4.27 in \cite{K02}) to conclude the convergence of the
above integrand to $\BE\exp\big[-\int f(s)\Pi_\rho(ds)\big]$.
Therefore the integral has this limit as well. It remains to prove that
$R(b)\to 0$ as $b\to\infty$. It is clearly sufficient to show that
$\BP_x(N^b_+(\tau(u)>0)=\BP_x(\tau(b)<\tau(u))\to 0$, where we recall
that $u>u_0$ and \eqref{4.51}. Let $\varepsilon>0$. As in the proof of 
Lemma \ref{l6} we choose a random walk with negative drift that is
dominating our process as long as it stays above $u$. We can then choose
$C>0$ large enough so that the maximum of this random walk is less than
$C$ with probability at least $1-\varepsilon$. Next we can choose $b>C$
large enough so that $\BP_x(X_{\tau^+_1(u)}>b-C)\le \varepsilon$.
This yields   $\BP_x(\tau(b)<\tau(u))\le 2\varepsilon$.\qed


\begin{remark}\label{rpos}\rm The positivity assumption in
Theorem \ref{geocompound} can be checked with the help of Lemma \ref{l2}.
To indicate how this can be done, we fix a $u\in \R$ satisfying
$\nu(u)>0$ (see Remark \ref{rem88}). 
Assume first that $q(u,t)\to\infty$ as $t\to\infty$. Since $\lambda$
is locally bounded we then have $\BP_u(\tau(b)<\infty)>0$.
Assume second that there are $\varepsilon,\delta>0$ such that
$\lambda(x)J(x,[\varepsilon,\infty))>0$ for all $x\ge u-\delta$. Due to the possibility
of many positive jumps in a small period of time we then have
$\BP_u(\tau^+_1(b)<\infty)>0$. If we now assume in addition that
$\lim_{t\to\infty}q(x,t)<u$ for all $x\ge u$,
then there is a positive probability for the process to drop
to level $b$ in a continuous way.
Hence we have again that $\BP_u(\tau(b)<\infty)>0$. 
\end{remark}

\section{Appendix}
\setcounter{equation}{0}

First we formulate some assumptions that will imply Assumptions \ref{a1}. 
Let us introduce the mean values
$$
m^-(x):=-\int_{-\infty}^0 z J(x,dz),\quad m^+(x):=\int_0^{\infty} z J(x,dz),
\quad x\in\R,
$$
and
$$
m(x):=m^+(x)-m^-(x)=\int z J(x,dz),\quad x\in\R.
$$

\begin{assumption}\label{a3}\rm
$m^-(x)+m^+(x)<\infty$ for all $x\in\R$ and 
$\lambda(x)(m^-(x)+m^+(x))$ is a locally bounded function on $\R$.
\end{assumption}



In the next assumption we use the convention $0/0:=0$.

\begin{assumption}\label{a5}\rm We have
\begin{align}
\label{C3}
\lim_{x\to-\infty}\frac{1}{m^+(x)}\int^{\infty}_{-x}(x+z)J(x,dz)
=\lim_{x\to\infty}\frac{1}{m^-(x)}\int_{-\infty}^{-x}(x+z)J(x,dz)=0.
\end{align}
\end{assumption}

Next we formulate a basic ergodicity assumption.

\begin{assumption}\label{a6}\rm There is an $\varepsilon>0$
such that
\begin{align}\label{C6}
\liminf_{x\to -\infty}(\mu(x)+\lambda(x)m^+(x)(1-\varepsilon)-\lambda(x)m^-(x))&>0,\\
\label{C7} 
\limsup_{x\to \infty}(\mu(x)+\lambda(x)m^+(x)-\lambda(x)m^-(x)(1-\varepsilon))&<0.
\end{align}
\end{assumption}

\begin{remark}\label{r7}\rm Assume that two of the limits
$\lim_{x\to-\infty}\mu(x)$, $\lim_{x\to-\infty}\lambda(x)m^-(x)$
and $\lim_{x\to-\infty}\lambda(x)m^+(x)$
exist and are finite, and make a similar assumption on the
corresponding limits as $x\to\infty$. Then Assumption \ref{a6} is equivalent to
\begin{align}\label{C61}
\liminf_{x\to -\infty}(\mu(x)+\lambda(x)m(x))>0>
\limsup_{x\to \infty}(\mu(x)+\lambda(x)m(x)).
\end{align}
For a constant (positive) $\mu$ and negative jumps
this is the well-known ergodicity condition
for the stress release model (see \cite{VereJones1988}, \cite{Zheng}, \cite{La04}).
\end{remark}

The next assumption is saying that all bounded sets are {\em small}
for the process (see \cite{MT1993}). Previous studies (see e.g.\
\cite{Zheng},\cite{La04},\cite{MT1993}) show that this is a rather
weak though sometimes tedious to check assumption. We will not
discuss it any further.

\begin{assumption}\label{a7}\rm For any bounded interval $I\subset \R$
there is a $t_0>0$ and a non-trivial  measure $\BQ$ on $\R$ such
that
$$
\BP_x(X_{t_0}\in\cdot)\ge \BQ(\cdot),\quad x\in I. 
$$
\end{assumption}

\begin{theorem}\label{terg}
If Assumptions \ref{a3}, \ref{a5}, \ref{a6} and \ref{a7} are satisfied,
then $\BP_x(T_\infty=\infty)=1$ for all $x\in\R$ and
$(X_t)$ has a unique invariant distribution $\pi$.
\end{theorem}
{\sc Proof:} We proceed similarly to \cite{La04}.
For any $m\ge 1$ the process $(X_{t\wedge\tau_m})$ is again Markov.
By \eqref{13} its generalized generator $\mathcal{A}_m$ (cf.\ \cite{MT1993})
is given by
\begin{align}\label{generator}
\mathcal{A}_mf(x)=\mu(x) f'(x)
+\lambda(x)\int (f(x+z)-f(x))J(x,dz),\quad |x|<m,
\end{align}
where $f$ satisfies the assumptions of Proposition \ref{p0}.
By Assumption \ref{a3} we can take $f(x):=|x|$ to obtain 
for $|x|<m$ that
\begin{align}\label{xc}\notag
\mathcal{A}_mf(x)=&\sgn(x)\mu(x) +\lambda(x)\int (|x+z|-|x|)J(x,dz)\\
=&\sgn(x)\mu(x) +\sgn(x)\lambda(x)m^+(x)-\sgn(x)\lambda(x)m^-(x)\\
\notag
&+2\bigg(\I\{x< 0\}\lambda(x)\int^\infty_{-x}(x+z)J(x,dz)
-\I\{x\ge 0\}\lambda(x)\int_{-\infty}^{-x}(x+z)J(x,dz)\bigg),
\end{align}
where $\sgn(x)\in\{-1,1\}$  is the sign of $x\in\R$, defined in a right-continuous way,
and where the second equality comes from 
$$
|x+z|-|x|=2(\I\{-z<x< 0\}-\I\{z<-x\le 0\})(x+z)+\sgn(x)z,\quad z\ne 0.
$$
We define 
$$
\varepsilon(x):=\I\{x< 0\} \frac{1}{2 m^+(x)}\int^{\infty}_{-x}(x+z)J(x,dz)
-\I\{x\ge 0\} \frac{1}{2 m^-(x)}\int_{-\infty}^{-x}(x+z)J(x,dz).
$$
Then $\varepsilon(x)\ge 0$ and from \eqref{C3} we have that
$\varepsilon(x)\to 0$ as $|x|\to\infty$. We can now rewrite \eqref{xc} as
\begin{align}\label{23}\notag
\mathcal{A}_mf(x)=&\sgn(x)\mu(x)
+\sgn(x)\lambda(x)m^+(x)(1-\I\{x< 0\}\varepsilon(x))\\
&-\sgn(x)\lambda(x)m^-(x)(1-\I\{x\ge 0\}\varepsilon(x)).
\end{align}
Using our assumptions in \eqref{23}, we easily get numbers
$\varepsilon> 0$, $x_0>0$, and $d\ge 0$ such that
\begin{align}\label{3214}
\mathcal{A}_mf(x)\le -\varepsilon+\I\{|x|\le x_0\}d,\quad |x|<m,m\in\N.
\end{align}
In particular, we can
apply Theorem 2.1 in \cite{MT1993} to conclude for any $x\in\R$ that
$\tau_m\to\infty$ $\BP_x$-almost surely as $m\to\infty$. This
proves the first assertion. We are then in a position to apply
Theorem 4.2 in \cite{MT1993} to complete the proof of the theorem.\qed

\begin{remark}\label{r41}\rm Under the conditions of Theorem \ref{terg},
the process $(X_t)$ is even {\em positive Harris recurrent}, see
\cite{MT1993}. Only a weak additional assumption is needed to 
obtain {\em Harris ergodicity}, i.e.\ the total variation convergence of
$\BP_x(X_t\in\cdot)$ to $\pi$ for any $x\in\R$.
By Theorem 6.1 in \cite{MT1993a}, one such assumption is {\em irreducibility}
of one {\em skeleton chain}.
\end{remark}


We next discuss Assumption \ref{a2}.
If $\lambda$ is a bounded function, then this assumption  is trivially satisfied.
If not, then we can impose the 
following slightly stronger version of Assumption \ref{a6} and
a weak positivity assumption on $m^-(x)+m^+(x)$.

\begin{assumption}\label{a8}\rm There is an $\varepsilon>0$
such that
\begin{align}\label{C62}
\liminf_{x\to -\infty}(\mu(x)+\lambda(x)m^+(x)(1-\varepsilon)-\lambda(x)m^-(x)(1+\varepsilon))&>0,\\
\label{C72} 
\limsup_{x\to \infty}(\mu(x)+\lambda(x)m^+(x)(1+\varepsilon)-\lambda(x)m^-(x)(1-\varepsilon))&<0.
\end{align}
\end{assumption}

\begin{theorem}\label{tergint}
If Assumptions \ref{a3}, \ref{a5}, \ref{a7} and \ref{a8} are satisfied
and, moreover,
\begin{align}\label{C8}
\liminf_{|x|\to\infty}(m^-(x)+m^+(x))>0,
\end{align}
then $\BP_x(T_\infty=\infty)=1$ for all $x\in\R$
and $(X_t)$ has a unique invariant distribution $\pi$ satisfying
$\int \lambda(x)\pi(dx)<\infty$.
\end{theorem}
{\sc Proof:} Using the assumptions in \eqref{23}, we can easily 
strengthen \eqref{3214} to
\begin{align}\label{3215}
\mathcal{A}_mf(x)\le -\max\{\varepsilon,\lambda(x)\}+\I\{|x|\le x_0\}d,\quad |x|<m,m\in\N.
\end{align}
Hence we can apply Theorem 4.2 in \cite{MT1993} to obtain that $\int\lambda(x)\pi(dx)<\infty$.
\qed

\begin{remark}\label{r14}\rm 
In the framework described in Remark \ref{r0}, Assumption \ref{a6}
can be reduced to \eqref{C7}. A similar remark applies to 
Assumptions \ref{a5} and \ref{a8}, and to \eqref{C8}.
\end{remark}

\vspace{0.3cm}
\noindent
{\bf Acknowledgement:} This research was supported by the ARC Centre
of Excellence for Mathematics and Statistics of Complex Systems.


\begin{thebibliography}{99}

\bibitem{Aldous89}
\newblock Aldous, D. (1989).
{\it Probability Approximations via the Poisson Clumping Heuristic.} 
Springer, New York.


\bibitem{As03}
Asmussen, S. (2003).
\newblock{\em Applied Probability and Queues}.
\newblock Second Edition, Springer, New York.


\bibitem{BB}
{\sc Baccelli, F. and Br\'emaud, P.} (1994).
{\it Elements of Queueing Theory.} Springer, Berlin.

\bibitem{BN72}
{\sc Bar-David, I. and Nemirovsky, A.} (1972).
Level crossings of nondifferentiable shot processes.
{\it IEEE Trans.\ Inform.\ Theory} {\bf 18}, 27--34.


\bibitem{Bo76}
{\sc Borovkov, A.A.} (1976). 
{\em Stochastic Processes in Queueing Theory.} Springer, New York.

\bibitem{BoNo01}
{\sc Borovkov, K.A., and Novikov, A.A.} (2001). 
On a piece-wise deterministic Markov process
model. {\em Stat. Probab. Letters,} {\bf 53}, 421--428.



\bibitem{BorovkovV}
{\sc Borovkov, K. and Vere-Jones, D.} (2000).
\newblock{Explicit formulae for stationary distributions of stress release processes.}
{\it Journal of Applied Probability} {\bf 37}, 2000, 315--321.


\bibitem{BrillPosner}
{\sc Brill, P.H.\ and Posner, M.J.M.} (1977)
Level crossings in point processes applied to queues: single server case.
{\it Operat.\ Res.} {\bf 25}, 662--673.


\bibitem{BrowneSigman}
{\sc Browne, S. and Sigman, K.} (1992).
Work-modulated queues with applications to storage processes.
{\it Journal of Applied Probability} {\bf 29}, 699--712.

\bibitem{Costa}
{\sc Costa, O.L.V.} (1990).
Stationary distributions for piecewise-deterministic Markov processes.
{\it Journal of Applied Probability} {\bf 27}, 60--73.


\bibitem{Davis}
{\sc Davis, M.H.A.} (1993). 
{\it Markov Models and Optimization.} Chapman and Hall, London.

\bibitem{Doshi92}
{\sc Doshi, B.T.} (1992). Level crossing analysis of queueing systems. 
in {\it Queueing and Related Models} Eds.\ Basawa, I.\ and Bhat U.N., 
Oxford Statsitical Science Series, Clarendon Press, Oxford.

\bibitem{FHR04}
{\sc Falk, M., H\"usler, J.\ and Reiss, R.-D.}
{\it Laws of Small Numbers: Extremes and Rare Events.}
Birkh\"auser, Basel.



\bibitem{HR}
{\sc Harrison, J.M. and Resnick, S.I.} (1976).
The stationary distribution and first exit probabilities of a storage process
with general release rule.
{\it Math. Oper. Res.} {\bf 1}, 347--358.


\bibitem{K02}
{\sc Kallenberg, O.} (2002).
\newblock{\it Foundations of Modern Probability}.
2nd Edition, Springer, New York.


\bibitem{La04}
{\sc Last, G.} (2004).
Ergodicity properties of stress release, repairable system
and workload models.
{\it Advances in Applied Probability} {\bf 36}, 471--498.


\bibitem{LastBrandt}
{\sc Last, G. and Brandt, A.} (1995). 
\newblock {\it Marked Point Processes on the Real Line.}
\newblock Probability and its Applications, Springer, New York.

\bibitem{LS98a}
{\sc Last, G. and Szekli, R.} (1998).
Stochastic comparison of repairable systems.
{\it Journal of Applied Probability} {\bf 35}, 348--370. 

\bibitem{Lead66}
{\sc Leadbetter, M.R.} (1966). 
On crossings of levels and curves by a wide class of stochastic processes.
{\it Annals of Mathematical Statistics} {\bf 37}, 260--267.

\bibitem{LeadHsing90}
{\sc Leadbetter, M. R. and Hsing, T.} (1990). 
Limit theorems for strongly mixing stationary random measures. 
{\it Stochastic Process. Appl.} {\bf 36}, 231--243.

\bibitem{LeadSpan02}
{\sc Leadbetter, M.R. and Spaniolo, G.V.} (2002).
On statistics at level crossings by a stationary process.
{\it Statistica Neerlandica} {\bf 56} (2), 152--164.

\bibitem{LinLeadRoo83}
{\sc Lindgren, G., Leadbetter, M.R., and Rootzen H.} (1983).
\newblock {\it Extremes and Related Properties of Stationary Sequences and Processes}. 
Springer, New York.


\bibitem{MT1993a}
{\sc Meyn, S.P. and Tweedie, R.L.} (1993).
Stability of Markovian processes II:
Continuous-time processes and sampled chains.
{\it Adv. Appl. Prob.} {\bf 25}, 487--517.


\bibitem{MT1993}
{\sc Meyn, S.P. and Tweedie, R.L.} (1993).
Stability of Markovian processes III:
Foster-Lyapunov criteria for continuous-time processes.
{\it Adv. Appl. Prob.} {\bf 25}, 518--548.

\bibitem{Rice}
{\sc Rice, S.O.} (1944).
Mathematical analysis of random noise.
{\it Bell System Tech. J.} {\bf 24}, 46--156.

\bibitem{VereJones1988}
{\sc Vere-Jones, D.} (1988).
On the variance properties of stress release models.
{\it Austral. J. Statist.} {\bf 30A}, 123--135. 


\bibitem{Zheng}
{\sc Zheng, X.} (1991).
Ergodic theorems for stress release processes.
{\it Stoch. Proces. Appl.} {\bf 37}, 239--258.



\end{thebibliography}
\end{document}